\documentclass[10pt,a4paper]{article}

\usepackage{graphicx,comment}
\usepackage{amsmath,amsthm}
\usepackage[numbers,sort&compress]{natbib}
\usepackage{microtype,multirow,bm}
\usepackage[bitstream-charter]{mathdesign}
\usepackage[final]{showkeys}
\usepackage[section]{placeins}
\usepackage{setspace}

\numberwithin{equation}{section}
\usepackage[hmargin=3cm,vmargin={3.5cm,4cm}]{geometry}
\allowdisplaybreaks

\newtheorem{Theorem}{Theorem}[section]
\newtheorem{Remark}{Remark}[section]

\title{Renewal processes based on generalized Mittag--Leffler waiting times}

\author{$\text{Dexter O. Cahoy}_1$, $\text{Federico Polito}_2$\\
	\footnotesize (1) -- Department of Mathematics and Statistics\\
    \footnotesize College of Engineering and Science, Louisiana Tech University, USA\\
    \footnotesize Tel: +1 318 257 3529, fax: +1 318 257 2182\\
	\footnotesize Email address: dcahoy@latech.edu\\
	\footnotesize (2) -- Department of Mathematics, University of Torino, Italy\\
	\footnotesize Tel: +39 011 6702937, fax: +39 011 6702878\\
	\footnotesize Email address: federico.polito@unito.it\\
	\footnotesize (Corresponding author)
	}

\begin{document}

	\maketitle
	
	\begin{abstract}

		The fractional Poisson process has recently attracted experts from several fields of study.
		Its natural generalization of the ordinary Poisson process made the model more appealing for real-world 
		applications. In this paper, we generalized the standard and fractional Poisson processes through the waiting time
		distribution, and showed their relations to an integral operator with a generalized Mittag--Leffler function
		in the kernel. The waiting times of the  proposed renewal processes have the generalized Mittag--Leffler
		and stretched-squashed Mittag--Leffler distributions.
		Note that the generalizations naturally provide greater flexibility in modeling real-life renewal processes.
		Algorithms to simulate sample paths and to estimate the model parameters are derived. Note also that these
		procedures are necessary to make these models more usable in practice.
		State probabilities and other qualitative or quantitative features of the models are also
		discussed.
		
		\vspace{.2cm}
		

	\end{abstract}
	
	\section{Introduction}

		The fractional Poisson process
		\citep{laskin,laskin2,sib,beg,beg2,cuw,mainardi,mainardi2,sibatov,meer,scalas,scalas2,scalas3}
		gained popularity in many areas of research as it naturally generalizes the standard or classical Poisson process.
		Recall that the inter-event time density function of the fractional Poisson process
		$N^\nu(t)$, $t \ge 0$, $\nu \in (0,1]$, was originally derived in
		\citet{ras00} (known to date) and has the following integral form:
		\begin{equation}
			\label{2e6}
			f^\nu(t)=\frac{1}{t}\int\limits_0^\infty
			e^{-x}\phi_\nu(\lambda t/x) \, \mathrm dx, \qquad \nu \in (0,1], \: t > 0, \lambda > 0,
		\end{equation}
		where
		\begin{align}
			\phi_\nu(\xi)=\frac{\sin(\nu\pi)}{\pi[\xi^\nu+\xi^{-\nu}+2\cos(\nu\pi)]}.
		\end{align}
		The preceding density function suggests that the tail distribution of the waiting time is of the form
		\begin{equation}
			\label{2e5}
			\Pr \, (T^\nu>t)= E_\nu(-\lambda t^\nu),
		\end{equation}
		where
		\begin{equation}
			\label{MLfunction}
			E_\beta(z)=\sum\limits_{n=0}^\infty\frac{z^n}{\Gamma(\beta n+1)},
			\qquad \: z \in \mathbb{C}, \: \beta \in \mathbb{C}, \Re (\beta)> 0,
		\end{equation}
		is the Mittag--Leffler function. Note that the Mittag--Leffler density has been widely used to describe
		distributions appearing in anomalous diffusion, finance and economics, transport of charge carriers in
		semiconductors, and light propagation through random media (see, e.g., \cite{uz99, psw05}).
		In view of equations \eqref{2e5} and \eqref{MLfunction}, the interarrival time density for the
		fractional Poisson process directly follows as
		\begin{equation}
			 \label{2e16}
			f^\nu(t)=\lambda t^{\nu-1} E_{\nu,\,\nu}(-\lambda t^{\nu}), \qquad t > 0,
		\end{equation}
		where
		\begin{align}
			E_{\beta,\gamma}(z) = \sum_{r=0}^{\infty}\frac{z^{r}}{\Gamma
			(\beta r + \gamma)}, \qquad z \in \mathbb{C}, \: \beta,\gamma \in \mathbb{C}, \: \Re(\beta)> 0
		\end{align}
		is the two-parameter Mittag--Leffler function. The $q$th fractional moment
		\citep{cap10} of the random interarrival  time is
		\begin{equation}
			\mathbb{E}\left[T^\nu\right]^q =\frac{\pi \Gamma (1 + q)}{\lambda^q \Gamma ( q / \nu )\sin ( \pi q/ \nu)
			\Gamma (1- q )}, \qquad 0< q < \nu.
		\end{equation}
		In addition, the above information automatically gives the probability density function
		\begin{equation}
			f_m^\nu(t)= \lambda^m  \frac{ t^{\nu m-1}}{(m-1)!}E_{\nu,\nu}^{(m-1)}\big(- \lambda t^{\nu} \big),
			\label{2e20}
		\end{equation}
		of the $m$-th arrival time  because its Laplace transform,
		\begin{align}
			\textsf{L} \big\lbrace f_m^\nu(t) \big\rbrace (s) =
			\int_0^\infty e^{-s t} f_m^\nu(t) \, \mathrm dt = \frac{\lambda^m}{(\lambda + s^\nu)^{m}},
		\end{align}
		where $ E_{\nu,\nu}^{(k)} (-\lambda t^\nu)$ is the $k$th derivative of $E_{\nu,\nu} (z)$
		evaluated at  $z =  -\lambda t^\nu$.  As $\nu \to 1$, the above distribution converges
		to the classical Erlang distribution.

		In another approach to the study the fractional Poisson process,
		\citet{laskin} used the fractional Kolmogorov--Feller-type differential equation system
		to characterize  the  one-dimensional state probability distributions as
		(see \citet[formula (25)]{laskin} and \citet[formula (2.5)]{beg2})
		\begin{align}
			\label{2e14}
			p_k^\nu(t) = \Pr \{ N^\nu(t) = k \} = \frac{(\lambda
			t^\nu)^{k}}{k!}\sum_{r=0}^{\infty}\frac{(r+k)!}{r!}\frac{(-\lambda
			t^{\nu})^r}{\Gamma (\nu (r+k) +1)}, \qquad
			k \ge 0, \: t \ge 0. 
		\end{align}
		One can also show \citep{laskin} that  the moment generating function
		(MGF) of the fractional Poisson process is
		\begin{align}
			M_\nu(s,t) & = E_\nu \left[ \lambda(e^{-s}-1)t^\nu \right] 
			= \sum_{r=0}^{\infty}\frac{\left[\lambda t^\nu\left(e^{-s}-1 \right)\right]^r}{\Gamma ( \nu  r+1)},
		\end{align}
		which permits calculation (see Table \ref{t1}) of the  moments.
		A summary of the characteristics of the classical and fractional Poisson processes is shown in Table \ref{t1}
		below.
		\begin{table}
			\label{t1}
	        \centering
			\begin{tabular}{c|cc}
				& Poisson process $(\nu=1)$ & Fractional Poisson Process $(\nu < 1)$\\
				\hline
				& & \\
				$\Pr(T^\nu>t)$& $e^{-\lambda t}$ & $E_{\nu}(-\lambda t^{\nu} )$\\
				& & \\
				$f^\nu (t)$&$ \lambda e^{-\lambda t} $& $\lambda t^{\nu-1} E_{\nu,\,\nu}(-\lambda t^{\nu})$\\
				& & \\
				$p_k^\nu(t)$& $\frac{(\lambda t)^{k}}{k!}e^{-\lambda t}$ & $\frac{(\lambda t^\nu)^{k}}
				{k!}\sum_{r=0}^{\infty}\frac{(r+k)!}{r!}\frac{(-\lambda t^{\nu})^r}{\Gamma (\nu (r+k) +1)}$\\
				& &\\
				Mean& $\lambda t$ &$\lambda t^{\nu}/\Gamma (\nu + 1)$\\
				& &\\
				Variance& $\lambda t$ & $\frac{\lambda t^\nu}{\Gamma(\nu+1)} + (\lambda t^\nu)^2
				\left[ \frac{1}{\nu\Gamma(2\nu)} - \frac{1}{\Gamma^2(\nu+1)} \right]$\\
				& &\\
				$k$th moment&$ (-1)^k \frac{\partial^k}{\partial s^k} \exp\left[ \lambda(e^{-s}-1)t\right]\big|_{s=0}$&$
				\left( -1\right)^k\frac{\partial^k}{\partial s^k} E_\nu\left[\lambda (e^{-s}-1)t^\nu\right]\big|_{s=0}$\\
				& & \\
			\end{tabular}
			\caption{\emph{Properties of fractional Poisson process compared with those of the standard Poisson
			process.}}
		\end{table}

		In this paper, we generalize the standard and fractional Poisson processes
		through their waiting time distributions.
		In particular, we propose two renewal processes that have waiting times that are generalized
		Mittag--Leffler and  stretched-squashed Mittag--Leffler distributed.
		These generalizations naturally provide more flexibility in capturing real-world renewal processes.
		Algorithms to simulate sample paths and estimate the model parameters are derived and tested.
		State probabilities and other qualitative or quantitative features of the models are also discussed.

		The rest of the paper is organized as follows. In Section \ref{00sec},
		a renewal process with generalized Mittag--Leffler distributed waiting times  is presented.
		Procedures to generate sample paths and to estimate
		parameters are also derived. In Section \ref{00anothersec}, another generalization based on
		stretching and squashing the Mittag--Leffler distributed inter-event times is developed.
		Methods  to simulate sample trajectories and to estimate parameters are also showcased.
		More discussions are provided in Section \ref{summary}.
		Finally, computational test results are shown in the appendix.

	\section{Generalization I}

		\label{00sec}	
\begin{figure}
			\centering
			\includegraphics[scale=.34]{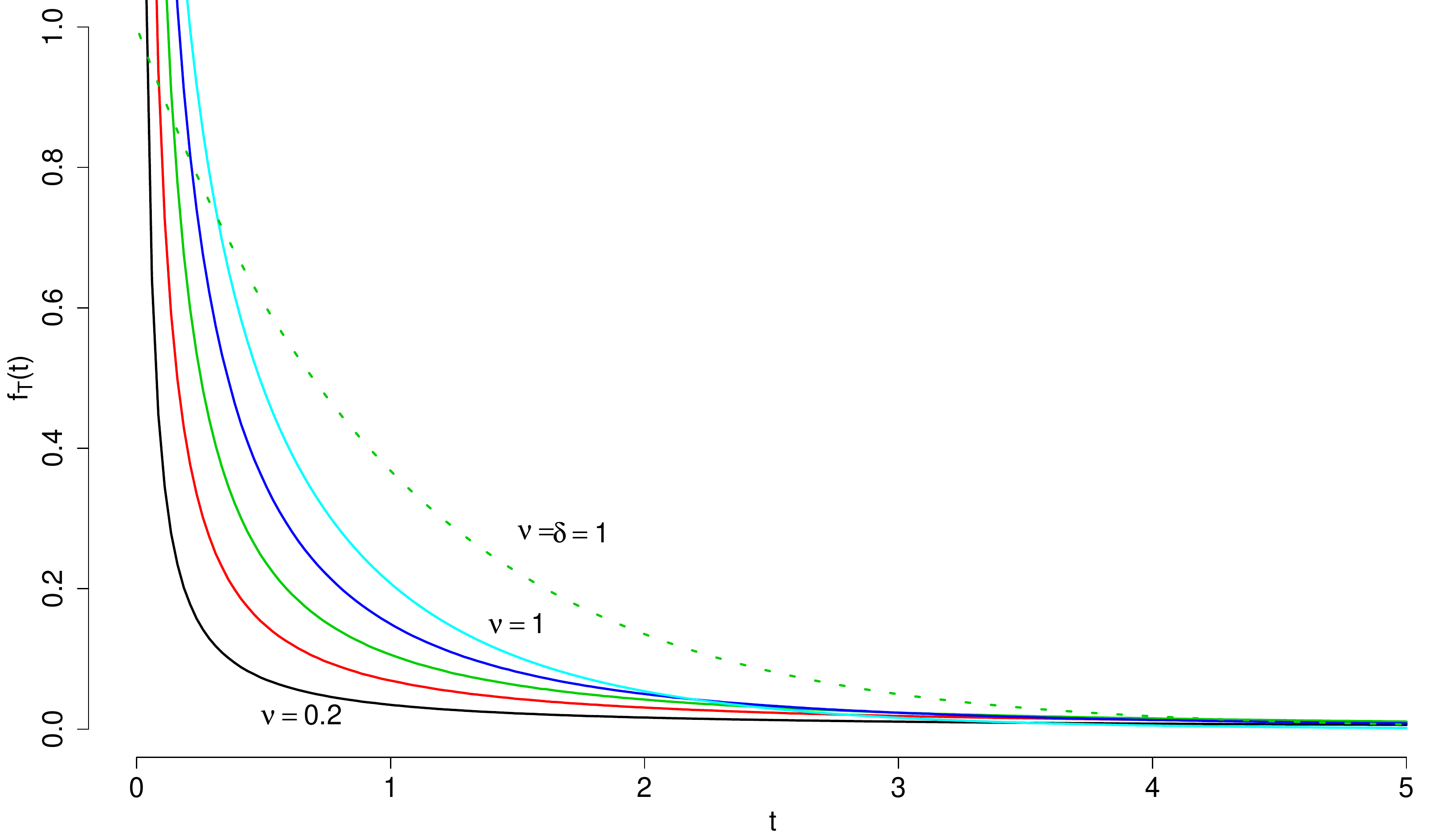}
			\includegraphics[scale=.34]{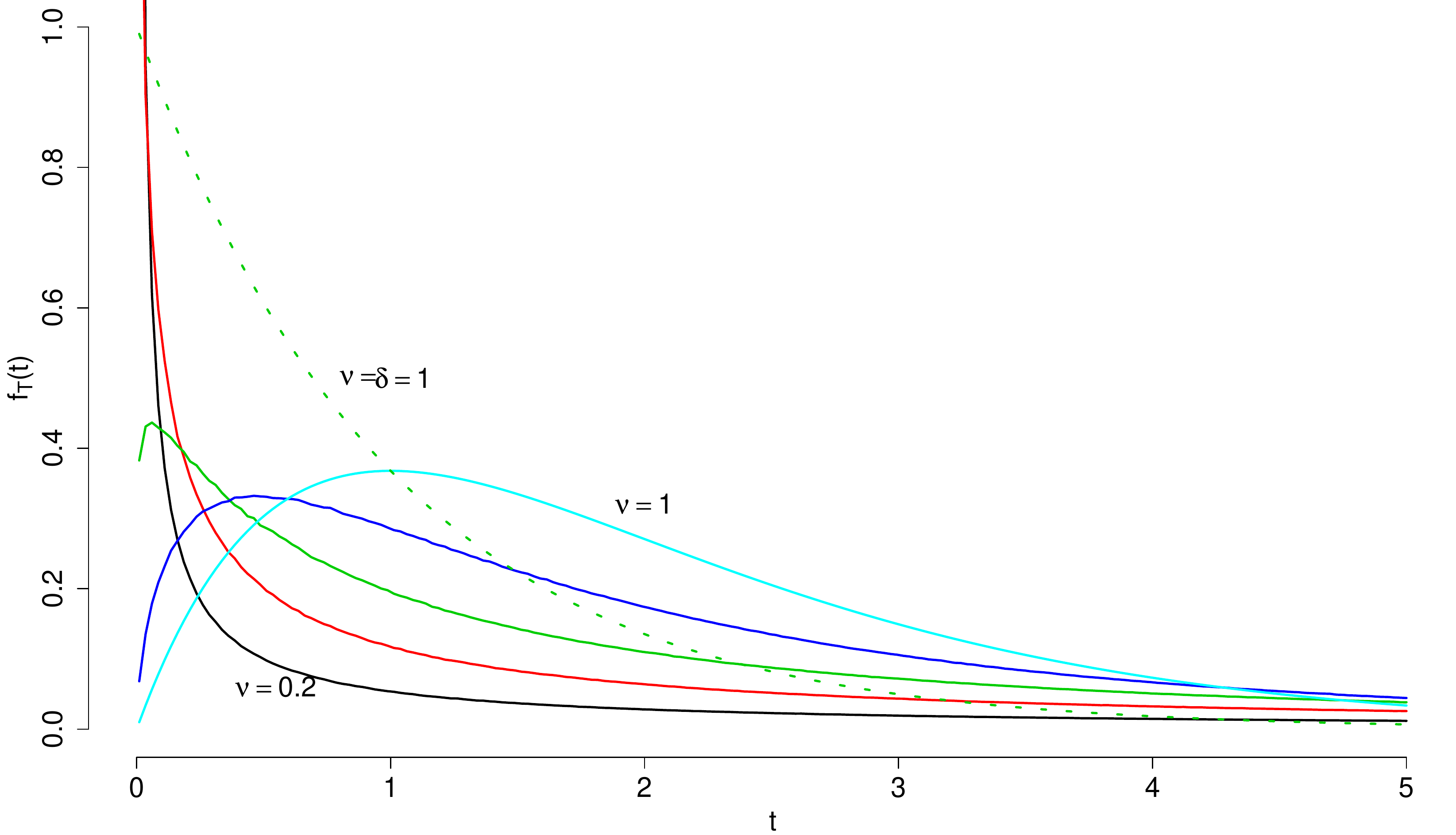}
			\caption{\label{f1}The generalized Mittag--Leffler density
			plots (see \eqref{00fu}) for parameter values $( \delta, \lambda )=( 0.5, 1)$ (top)
			and $(\delta, \lambda )=(2, 1)$ (bottom) where $\nu$  goes from 0.2 to 1 with step size of 0.2.}
		\end{figure}

		We consider the generalized Mittag--Leffler distribution (see e.g.\ \citet{pil})
		built from the generalized Mittag--Leffler function \citep{saigo,prab}. Let $T^{\nu,\delta}$ be a
		generalized Mittag--Leffler distributed random
		variable. Then the probability density function is
		\begin{align}
			\label{00fu}
			f^{\nu,\delta}(t) = \lambda^\delta t^{\delta \nu -1} E_{\nu, \delta \nu}^\delta (-\lambda t^\nu),
			\qquad t>0, \: \lambda > 0, \: \nu \in (0, 1], \: \delta \in \mathbb{R},
		\end{align}
		where
		\begin{align}
			E_{\beta,\gamma}^\xi(z) = \sum_{r=0}^\infty \frac{(\xi)_r}{r!\Gamma(\beta r+\gamma)} z^r,
			\quad \beta,\gamma,\xi, z \in \mathbb{C}, \: \Re(\beta)>0
		\end{align}
		is the generalized Mittag--Leffler function (see Figure \ref{f1}). The Pochhammer symbol $(\xi)_r$ can be
		written also as $(\xi)_r = \xi(\xi+1)\dots (\xi+r-1)$,
		$\xi \ne 0$. When $\delta \nu < 1$ the function \eqref{00fu} has an asymptote at $t=0$, while
		in the particular case $\delta \nu = 1$
		\begin{align}
			\left. f^{\nu,\delta}(t) \right|_{t=0} = \left. \lambda^{1/\nu}
			E_{\nu,1}^{1/\nu}(-\lambda t^\nu) \right|_{t=0} = \lambda^{1/\nu}.
		\end{align}
		The Laplace transform of \eqref{00fu} reads
		\begin{align}
			\textsf{L} \big\lbrace \lambda^\delta t^{\delta\nu-1} E_{\nu,\delta\nu}^\delta (-\lambda t^\nu) \big\rbrace
			(s)= \frac{\lambda^\delta}{(s^\nu+\lambda)^\delta}
		\end{align}
		(see \citet{mathai}, formula (2.3.24), page 95). Below are the plots of the generalized Mittag--Leffler densities.
		
		Now, consider $m$ i.i.d.\ random waiting times $\mathcal{T}_i$, $i=1,\dots, m$
		of a renewal point process, here denoted as $N^{\nu,\delta}(t)$, $t \ge 0$,
		which are distributed as in \eqref{00fu}. 	Furthermore,
		denote $T_m^{\nu,\delta} = \mathcal{T}_1+\dots +\mathcal{T}_m$  as the waiting time of
		the $m$th renewal event and $T_0^{\nu,\delta} = 0$.  Then
		\begin{align}
			\mathbb{E} e^{-s T_m^{\nu,\delta}} = \left( \mathbb{E}e^{-s\mathcal{T}_i} \right)^m 
			= \frac{\lambda^{\delta m}}{(s^\nu+\lambda)^{\delta m}}.
		\end{align}
		By formula (2.3.24) of \citet{mathai} we have
		\begin{align}
			\Pr \{ T_m^{\nu,\delta} \in \mathrm dt \} / \mathrm dt = \lambda^{\delta m} t^{\nu \delta m-1}
			E_{\nu,\nu\delta m}^{\delta m} (-\lambda t^\nu), \qquad t > 0, \: \lambda > 0.
		\end{align}
		It is  rather immediate now to obtain the state probabilities
		$p_k^{\nu,\delta}(t) = \Pr \{ N^{\nu,\delta}(t) = k \}$, $k \ge 0$ because 
		\begin{align}
			\label{00state}
			& \int_0^\infty e^{-st} \Pr \{ N^{\nu,\delta}(t) = k \} \, \mathrm dt \\
			& =	\int_0^\infty e^{-st} \left( \Pr\{ T_k^{\nu,\delta}<t \} -
			\Pr \{ T_{k+1}^{\nu,\delta}<t \} \right) \mathrm dt \notag \\
			& = \int_0^\infty e^{-st} \left\{ \int_0^t \Pr\{ T_k^{\nu,\delta} \in \mathrm dy \}
			- \int_0^t \Pr\{ T_{k+1}^{\nu,\delta} \in \mathrm dy \} \right\} \, \mathrm dt \notag \\
			& = \int_0^\infty \Pr\{ T_k^{\nu,\delta} \in \mathrm dy \} \int_y^\infty e^{-st} \mathrm dt
			- \int_0^\infty \Pr\{ T_{k+1}^{\nu,\delta} \in \mathrm dy \} \int_y^\infty e^{-st} \mathrm dt \notag \\
			& = s^{-1} \left[ \int_0^\infty e^{-sy} \Pr \{ T_k^{\nu,\delta} \in \mathrm dy \} - \int_0^\infty e^{-sy} \Pr
			\{ T_{k+1}^{\nu,\delta} \in \mathrm dy \} \right] \notag \\
			& = s^{-1} \left( \frac{\lambda^{\delta k}}{(s^\nu+\lambda)^{\delta k}}
			-\frac{\lambda^{\delta(k+1)}}{(s^\nu+\lambda)^{\delta(k+1)}} \right) \notag, \qquad k \ge 0.
		\end{align}
		Inverting the preceding Laplace transform we readily arrive at
		\begin{align}
			\label{001}
			p_k^{\nu,\delta}(t) =
			\lambda^{\delta k} t^{\nu\delta k} E_{\nu,\nu\delta k+1}^{\delta k} (-\lambda t^\nu)
			-\lambda^{\delta (k+1)} t^{\nu\delta(k+1)} E_{\nu,\nu\delta(k+1)
			+1}^{\delta(k+1)} (-\lambda t^\nu), \qquad k \ge 0.
		\end{align}
		When $(0)_r=0$, $r \in \mathbb{N}$, $(0)_0=1$, and $k=0$, equation \eqref{001}  becomes
		\begin{align}
			p_0^{\nu,\delta}(t) = 1-\lambda^\delta t^{\nu\delta} E_{\nu,\nu\delta+1}^\delta (-\lambda t^\nu).
		\end{align}
		Clearly,
		\begin{align}
			\label{40p}
			p_k^{\nu,\delta}(0) = \Pr \{ N^{\nu,\delta}(0) = k\} =
			\begin{cases}
				1, & k = 0, \\
				0, & k \ge 1.
			\end{cases}
		\end{align}
		
		\begin{Theorem}
			The state probabilities $p_k^{\nu,\delta}(t) = \Pr \{ N^{\nu,\delta}(t) = k \}$,
			$k \ge 1$, satisfy the convolution-type Volterra
			equation of the first kind
			\begin{align}
				\label{00hat}
				p_k^{\nu,\delta}(t) = \lambda^\delta \int_0^t (t-w)^{\nu\delta -1}
				E_{\nu,\nu\delta}^\delta \left(-\lambda(t-w)^\nu \right)
				p_{k-1}^{\nu,\delta}(w) \, \mathrm dw.
			\end{align}
			\begin{proof}
				We start by rewriting \eqref{00state} by means of the following well-known relation
				(see e.g.\ \citet{sax}, formula (11.7), page 17):
				\begin{align}
					\label{01for}
					\int_0^x (x-t)^{\beta-1} E_{\alpha,\beta}^\gamma [a(x-t)^\alpha] t^{\nu-1}
					E_{\alpha,\nu}^\sigma (at^\alpha) \mathrm dt = x^{\beta + \nu-1} E_{\alpha, \beta+\nu}^{\gamma+\sigma}
					(ax^\alpha),
				\end{align}
				where $\alpha,\beta,\gamma,a,\nu,\sigma \in \mathbb{C}$, and $\Re(\alpha)>0$, $\Re(\beta)>0$,
				$\Re(\gamma)>0$, $\Re(\nu)>0$, $\Re(\sigma)>0$.
				Then
				\begin{align}
					p_k^{\nu,\delta} (t)
					= {} & \lambda^{\delta k} t^{\nu\delta k} E_{\nu,\nu\delta k +1}^{\delta k}(-\lambda t^\nu)
					- \lambda^{\delta(k+1)} t^{\nu\delta(k+1)} E_{\nu,\nu\delta(k+1)+1}^{\delta(k+1)}(-\lambda t^\nu) \\
					= {} & \lambda^{\delta k} \int_0^t (t-s)^{\nu\delta (k-1)} E_{\nu,\nu\delta(k-1)+1}^{\delta(k-1)}
					\left[ -\lambda (t-s)^\nu \right] s^{\nu\delta -1} E_{\nu,\nu\delta}^\delta (-\lambda s^\nu)
					\, \mathrm ds
					\notag \\
					& - \lambda^{\delta(k+1)} \int_0^t (t-s)^{\nu\delta k} E_{\nu,\nu\delta k +1}^{\delta k}
					\left[ -\lambda(t-s)^\nu \right] s^{\nu\delta -1}
					E_{\nu,\nu\delta}^\delta (-\lambda s^\nu) \, \mathrm ds \notag \\
					& \!\!\!\!\!\!\!\!\!\!\!\!\!\! \overset{w=t-s}{=}  \lambda^\delta \int_0^t (t-w)^{\nu\delta -1}
					E_{\nu,\nu\delta}^\delta \left[ -\lambda(t-w)^\nu \right]
					\lambda^{\delta(k-1)} w^{\nu\delta(k-1)} E_{\nu,\nu\delta(k-1)+1}^{\delta(k-1)}
					(-\lambda w^\nu) \, \mathrm dw \notag \\
					& - \lambda^\delta \int_0^t (t-w)^{\nu\delta -1} E_{\nu,\nu\delta}^\delta \left[
					-\lambda(t-w)^\nu \right]
					\lambda^{\delta k} w^{\nu\delta k} E_{\nu,\nu\delta k+1}^{\delta k} (-\lambda w^\nu)
					\, \mathrm dw \notag \\
					= {} & \lambda^\delta \int_0^t (t-w)^{\nu\delta -1}
					E_{\nu,\nu\delta}^\delta \left[ -\lambda (t-w)^\nu \right] \notag \\
					& \times \left[ \lambda^{\delta(k-1)} w^{\nu\delta(k-1)} E_{\nu,\nu\delta(k-1)+1}^{\delta(k-1)}
					- \lambda^{\delta k} w^{\nu\delta k} E_{\nu,\nu\delta k +1}^{\delta k}
					(-\lambda w^\nu) \right] \mathrm dw \notag \\
					= {} & \lambda^\delta \int_0^t (t-w)^{\nu\delta -1}
					E_{\nu,\nu\delta}^\delta \left[ (-\lambda(t-w)^\nu) \right]
					p_{k-1}^{\nu,\delta}(w) \, \mathrm dw. \notag
				\end{align}
			\end{proof}
		\end{Theorem}
		
		\begin{Remark}		
			Result \eqref{00hat} can also be conveniently expressed  by means of the Prabhakar operator \citep{prab},
			defined as
			\begin{align}
				\left( \bm{\mathrm E}^\gamma_{\rho,\mu,\omega;a+} \phi \right)(x) =
				\int_a^x (x-y)^{\mu-1} E_{\rho,\mu}^\gamma
				\left(\omega(x-y)^\rho\right) \phi(y) \, \mathrm dy, \qquad x>a, \: \rho,\mu,\gamma \in \mathbb{C}, \:
				\Re (\rho), \Re(\mu) >0,
			\end{align}
			which is a generalization of the Riemann--Liouville fractional integral of  $\phi (x)$.
			Therefore, we obtain
			\begin{align}
				\label{arancione}
				p_k^{\nu,\delta}(t)
				= \lambda^\delta \left( \bm{\mathrm E}_{\nu,\nu\delta,-\lambda;0+}^\delta \,
				p_{k-1}^{\nu,\delta} \right) (t).
			\end{align}
		\end{Remark}
		
		\begin{Remark}
			When $\delta = 1$ equation \eqref{00hat} clearly reduces to
			\begin{align}
				p_k^\nu(t) = \lambda \int_0^t (t-w)^{\nu-1} E_{\nu,\nu}(-\lambda(t-w)^\nu)
				p_{k-1}^\nu(w) \, \mathrm dw.
			\end{align}
			We now check that the state probabilities $p_k^\nu(t)$ of the fractional Poisson process
			$N^\nu(t)$, $t \ge 0$ (see e.g.\ \citet{beg2}) satisfy the above integral equation. By recalling that
			\begin{align}
				p_k^\nu(t) = (\lambda t^\nu)^k E_{\nu,\nu k+1}^{k+1} (-\lambda t^\nu), \qquad k \ge 0, \:
				t \ge 0,
			\end{align}
			we can write
			\begin{align}
				p_k^\nu(t) = {} & \lambda \int_0^t (t-w)^{\nu-1} E_{\nu,\nu}(-\lambda(t-w)^\nu)
				\, p_{k-1}^\nu \, \mathrm dw \\
				= {} & \lambda \int_0^t (t-w)^{\nu-1} E_{\nu,\nu}(-\lambda(t-w)^\nu)
				\, (\lambda w)^{k-1} E_{\nu,\nu(k-1)+1}(-\lambda w^\nu) \, \mathrm dw \notag \\
				& \hspace{-.76cm} \overset{\text{by \eqref{01for}}}{=}
				(\lambda t^\nu)^k E_{\nu,\nu k+1}^{k+1} (-\lambda t^\nu) \notag.
			\end{align}
						
			Notice also that for $\nu=\delta=1$ (classical case),  equation 	\eqref{00hat} reduces to
			\begin{align}
				\label{00eq}
				p_k(t) = \lambda \int_0^t e^{-\lambda(t-w)}
				p_{k-1}(w) \, \mathrm dw,
			\end{align}
			where $p_k(t)$, $k \ge 1$, $t \ge 0$, are the state probabilities of a
			homogeneous Poisson process $N(t)$, $t \ge 0$. Equation \eqref{00eq} is easily solvable and the solution
			reads
			\begin{align}
				p_{k-1}(t) = \lambda^{-1} \frac{\mathrm d}{\mathrm dt} p_k(t) +
				p_k(t), \qquad k \ge 1.
			\end{align}
			Finally we note that the above equation is clearly the difference-differential equation
			governing the state probabilities of a homogeneous Poisson process.
		\end{Remark}
		
		\begin{Theorem}
			\label{teio1}
			The state probabilities $p_k^{\nu,\delta}(t) = \text{Pr} \{ N^{\nu,\delta}(t) = k \}$, $k \ge 0$, $t \ge 0$,
			satisfy the equations
			\begin{align}
				\label{arancione1}
				\frac{\mathrm d^{\nu \delta +\theta}}{\mathrm dt^{\nu\delta+\theta}}
				\left( \bm{\mathrm E}_{\nu,\theta,-\lambda;0+}^{-\delta} \, p_k^{\nu,\delta} \right) (t)
				= \lambda^\delta p_{k-1}^{\nu,\delta}(t)
				+\delta_{k,0} \left( t^{-\nu\delta} E_{\nu,1-\nu\delta}^{-\delta}(-\lambda t^\nu) -\lambda^\delta \right),
			\end{align}
			for any $\theta \in \mathbb{C}$, $\Re (\theta)>0$, where $\delta_{k,0}$ is the Kronecker's delta and
			where the operator $\frac{\mathrm d^{\nu \delta +\theta}}{\mathrm dt^{\nu\delta+\theta}}$
			is the Riemann--Liouville fractional derivative of order $\nu\delta+\theta$.
			
			\begin{proof}
				We start by considering $k\ge 1$. Applying the  the operator
				$\bm{\mathrm E}_{\nu,\theta,-\lambda;0+}^{-\delta}$  to equation \eqref{arancione}, we obtain
				\begin{align} 
					& \left( \bm{\mathrm E}_{\nu,\theta,-\lambda;0+}^{-\delta} \, p_k^{\nu,\delta}
					\right) (t) = \lambda^\delta
					\left( \bm{\mathrm E}_{\nu,\theta,-\lambda;0+}^{-\delta}
					\left( \bm{\mathrm E}_{\nu,\nu\delta,-\lambda;0+}^\delta \, p_{k-1}^{\nu,\delta}
					\right) \right) (t) \\
					& \Leftrightarrow \quad
					\left( \bm{\mathrm E}_{\nu,\theta,-\lambda;0+}^{-\delta} \, p_k^{\nu,\delta} \right) (t) =
					\lambda^\delta J^{\nu\delta+\theta}_{t,0+} p_{k-1}^{\nu,\delta}(t), \notag
				\end{align}
				where $J^{\nu\delta+\theta}_{t,0+}$ is the Riemann--Liouville fractional integral operator.
				By recalling that the Riemann--Liouville fractional derivative is the left
				inverse operator of the Riemann--Liouville fractional integral (see e.g.\ \citet{Diethelm},
				Theorem 2.14, page 30), we readily arrive at the claimed result.
				For $k=0$,  it is sufficient to show that
				\begin{align}
					\vspace{-1cm} \frac{\mathrm d^{\nu \delta +\theta}}{\mathrm dt^{\nu\delta+\theta}} &
					\left( \bm{\mathrm E}_{\nu,\theta,-\lambda;0+}^{-\delta} \, p_0^{\nu,\delta} \right) (t) \\
					= {} & \frac{\mathrm d^{\nu \delta +\theta}}{\mathrm dt^{\nu\delta+\theta}}
					\left( \int_0^t (t-y)^{\theta-1} E_{\nu,\theta}^{-\delta} \left[ -\lambda(t-y)^\nu \right]
					p_0^{\nu,\delta}(y) \, \mathrm dy \right) \notag \\
					= {} & \frac{\mathrm d^{\nu \delta +\theta}}{\mathrm dt^{\nu\delta+\theta}}
					\left( \int_0^t (t-y)^{\theta-1} E_{\nu,\theta}^{-\delta} \left[ -\lambda(t-y)^\nu \right]
					\, \mathrm dy \right. \notag \\
					& \left. - \lambda^\delta \int_0^t (t-y)^{\theta-1} E_{\nu,\theta}^{-\delta}
					\left[ -\lambda (t-y)^\nu \right] y^{\nu\delta} E_{\nu,\nu\delta+1}^\delta (-\lambda y^\nu)
					\, \mathrm dy \right) \notag \\
					= {} & \frac{\mathrm d^{\nu \delta +\theta}}{\mathrm dt^{\nu\delta+\theta}}
					\left( t^\theta E_{\nu,\theta+1}^{-\delta} (-\lambda t^\nu)
					-\lambda^\delta t^{\nu\delta+\theta} E_{\nu,\nu\delta+\theta+1}^0 (-\lambda t^\nu) \right) \notag \\
					= {} & \frac{\mathrm d^{\nu \delta +\theta}}{\mathrm dt^{\nu\delta+\theta}}
					\left( t^\theta E_{\nu,\theta+1}^{-\delta} (-\lambda t^\nu)
					-\lambda^\delta \frac{t^{\nu\delta+\theta}}{\Gamma(\nu\delta +\theta+1)} \right) \notag \\
					= {} & t^{-\nu\delta} E_{\nu,1-\nu\delta}^{-\delta} (-\lambda t^\nu) -\lambda^\delta. \notag
				\end{align}
			\end{proof}
		\end{Theorem}
		For more information on the inverse operator appearing in \eqref{arancione1}, the reader can consult
		\citet{saigo}, Section 6.
		
		\begin{Remark}
			When $\delta=1$ and $k\ge 1$, equation \eqref{arancione1} can be written as
			\begin{align}
				p_{k-1}^\nu (t) & = \lambda^{-1} \frac{\mathrm d^{\nu+\theta}}{\mathrm d t^{\nu+\theta}}
				\int_0^t (t-w)^{\theta-1} E_{\nu,\theta}^{-1} (-\lambda(t-w)^\nu) p_k^\nu(t) \, \mathrm dw \\
				& = \lambda^{-1} \frac{\mathrm d^{\nu+\theta}}{\mathrm d t^{\nu+\theta}}
				\int_0^t (t-w)^{\theta-1} \left[ \frac{1}{\Gamma(\theta)} + \lambda(t-w)^\nu \frac{1}{\Gamma(\nu+\theta)}
				\right] p_k^\nu(w) \, \mathrm dw \notag \\
				& = \lambda^{-1} \frac{\mathrm d^\nu}{\mathrm dt^\nu} \frac{\mathrm d^\theta}{\mathrm dt^\theta}
				\frac{1}{\Gamma(\theta)} \int_0^t (t-w)^{\theta-1} p_k^\nu(w) \, \mathrm dw
				+ \frac{\mathrm d^{\nu+\theta}}{\mathrm dt^{\nu+\theta}} \frac{1}{\Gamma(\nu+\theta)} \int_0^t
				(t-w)^{\nu+\theta-1} p_k^\nu(w) \, \mathrm dw \notag,
			\end{align}
			while, for $k=0$, and considering that
			\begin{align}
				t^{-\nu} E_{\nu,1-\nu}^{-1} (-\lambda t^\nu) -\lambda & = t^{-\nu} \left( \frac{1}{\Gamma(1-\nu)}
				+\lambda t^\nu \right) -\lambda
				= \frac{t^{-\nu}}{\Gamma(1-\nu)},
			\end{align}
			we have
			\begin{align}
				\frac{t^{-\nu}}{\Gamma(1-\nu)} & = \frac{\mathrm d^{\nu+\theta}}{\mathrm d t^{\nu+\theta}}
				\int_0^t (t-w)^{\theta-1} E_{\nu,\theta}^{-1} (-\lambda(t-w)^\nu) \, p_0^\nu(t) \, \mathrm dw \\
				& = \frac{\mathrm d^{\nu+\theta}}{\mathrm d t^{\nu+\theta}}
				\int_0^t (t-w)^{\theta-1} \left[ \frac{1}{\Gamma(\theta)} + \lambda(t-w)^\nu \frac{1}{\Gamma(\nu+\theta)}
				\right] p_0^\nu(w) \, \mathrm dw \notag \\
				& = \frac{\mathrm d^\nu}{\mathrm dt^\nu} \frac{\mathrm d^\theta}{\mathrm dt^\theta}
				\frac{1}{\Gamma(\theta)} \int_0^t (t-w)^{\theta-1} p_0^\nu(w) \, \mathrm dw
				+ \lambda \frac{\mathrm d^{\nu+\theta}}{\mathrm dt^{\nu+\theta}} \frac{1}{\Gamma(\nu+\theta)} \int_0^t
				(t-w)^{\nu+\theta-1} p_0^\nu(w) \, \mathrm dw \notag.
			\end{align}
			Hence, we retrieve the fractional difference-differential equations governing the state probabilities
			of a fractional Poisson process \citep{laskin}:
			\begin{align}
				\label{pejo}
				\frac{\mathrm d^\nu}{\mathrm dt^\nu} p_k^\nu(t)
				= -\lambda p_k^\nu(t) + \lambda p_{k-1}^\nu(t) + \delta_{k,0}
				\frac{t^{-\nu}}{\Gamma(1-\nu)}, \qquad k \ge 0.
			\end{align}
		\end{Remark}
		
		From equation \eqref{arancione1} we can easily arrive at the following partial differential equation
		for the probability generating function $\mathcal{G}_{\nu,\delta}(u,t) = \sum_{k=0}^\infty u^k
		p_k^{\nu,\delta}(t)$.
		\begin{align}
			\label{12gen}
			\frac{\partial^{\nu\delta+\theta}}{\partial t^{\nu\delta+\theta}}
			\left( \bm{\mathrm E}_{\nu,\theta,-\lambda;0+}^{-\delta} \mathcal{G}_{\nu,\delta}(u,\cdot) \right) (t)
			= \lambda^\delta u \, \mathcal{G}_{\nu,\delta}(u,t)
			+ t^{-\nu\delta} E_{\nu,1-\nu\delta}^{-\delta} (-\lambda t^\nu) -\lambda^\delta.
		\end{align}
		From the above equation and by recalling the formula
		$\frac{\partial}{\partial u} \mathcal{G}_{\nu,\delta}(u,t)|_{u=1} = \mathbb{E} \,
		N^{\nu,\delta}(t)$,
		it is now immediate to derive the differential equation involving the mean value as
		\begin{align}
			\label{12mean}
			\frac{\mathrm d^{\nu\delta+\theta}}{\mathrm d t^{\nu \delta + \theta}}
			\left( \bm{\mathrm E}_{\nu,\delta,-\lambda;0+}^{-\delta} \mathbb{E}\, N^{\nu,\delta}(\cdot) \right) (t)
			= \lambda^\delta \left( 1+\mathbb{E} \, N^{\nu,\delta}(t)\right).
		\end{align}
		Observe that equations \eqref{12gen} and \eqref{12mean} reduce to the corresponding equations in the pure
		fractional case when $\delta=1$  (see \citet[formula (22)]{laskin}
		for the differential equation involving the probability generating function).
		For the fractional Poisson process  $N^\nu(t)$, $t \ge 0$  ($\delta= 1$),  equation \eqref{12mean}  becomes
		\begin{align}
			\label{30water}
			\frac{\mathrm d^{\nu+\theta}}{\mathrm d t^{\nu+\theta}} & \left( \bm{\mathrm{E}}_{\nu,\theta,-\lambda;0+}^{-1}
			\mathbb{E} \, N^\nu(\cdot) \right)(t) = \lambda + \lambda \mathbb{E} \, N^\nu(t) \\
			\Leftrightarrow {} \qquad & \frac{\mathrm d^{\nu+\theta}}{\mathrm d t^{\nu+\theta}}
			\int_0^t (t-y)^{\theta-1} E_{\nu,\theta}^{-1} [-\lambda (t-y)^\nu] \, \mathbb{E} \, N^\nu(t) \, \mathrm dy
			= \lambda + \lambda \mathbb{E} \, N^\nu(t) \notag \\
			\Leftrightarrow {} \qquad & \frac{\mathrm d^\nu}{\mathrm d t^\nu} \frac{\mathrm d^\theta}{\mathrm dt^\theta}
			\frac{1}{\Gamma(\theta)} \int_0^t (t-y)^{\theta-1} \mathbb{E} \, N^\nu(y) \, \mathrm dy \notag \\
			& + \lambda \frac{\mathrm d^{\nu+\theta}}{\mathrm dt^{\nu+\theta}} \frac{1}{\Gamma(\nu+\theta)}
			\int_0^t (t-y)^{\nu+\theta-1} \mathbb{E}\, N^\nu(y)\, \mathrm dy
			= \lambda + \lambda \mathbb{E} \, N^\nu(t) \notag \\
			\Leftrightarrow {} \qquad & \frac{\mathrm d^\nu}{\mathrm dt^\nu} \mathbb{E} \, N^\nu(t)
			= \lambda, \notag
		\end{align}
		with $\mathbb{E} \, N^\nu(0) = 0$ and considering that the second step is justified by the
		semigroup property of the Riemann--Liouville fractional derivative (see \citet[Theorem 2.2, page 14]{Diethelm}).
		The solution to \eqref{30water} is well-known and reads \citep[formula (26)]{laskin}
		\begin{align}
			\label{meanpure}
			\mathbb{E}\, N^\nu(t) = \frac{\lambda t^\nu}{\Gamma(\nu+1)}, \qquad \nu \in (0,1] \: t \ge 0.
		\end{align}
		
		The following theorem derives the mean value of the process $N^{\nu,\delta}(t)$, $t \ge 0$.
		\begin{Theorem}
			Let $\nu \in (0,1]$, $\delta \in \mathbb{C}$, $\theta \in \mathbb{C}$, $\Re (\theta)>0$.
			The solution to
				\begin{align}
					\label{30mean}
					\begin{cases}
						\frac{\mathrm d^{\nu\delta+\theta}}{\mathrm d t^{\nu \delta + \theta}}
						\left( \bm{\mathrm E}_{\nu,\theta,-\lambda;0+}^{-\delta} \mathbb{E}\, N^{\nu,\delta}
						(\cdot) \right) (t)
						= \lambda^\delta \left( 1+\mathbb{E} \, N^{\nu,\delta}(t)\right), \\
						\left[ \frac{\mathrm d^{\nu\delta +\theta -k-1}}{\mathrm dt^{\nu\delta +\theta -k-1}}
						\left( \bm{\mathrm E}_{\nu,\theta,-\lambda;0+}^{-\delta}
						\mathbb{E} \, N^{\nu,\delta} (\cdot) \right) (t) \right]_{t \to 0} =0,
						\qquad \forall k= 0, \dots, n-1,
						\: n-1 \le \Re (\nu\delta+\theta) < n,
					\end{cases}
				\end{align}
				reads
				\begin{align}
					\label{30meanres}
					\mathbb{E} \, N^{\nu,\delta}(t) = \sum_{r=0}^\infty \lambda^{\delta(r+1)} t^{\nu\delta (r+1)}
					E_{\nu,\nu\delta(r+1)+1}^{\delta(r+1)}(-\lambda t^\nu).
				\end{align}
				
				\begin{proof}
					We start by taking the Laplace transform of \eqref{30mean}, obtaining
					\begin{align}
						\label{30mouse}
						& \int_0^\infty e^{-st}
						\frac{\mathrm d^{\nu\delta+\theta}}{\mathrm d t^{\nu \delta + \theta}}
						\left( \bm{\mathrm E}_{\nu,\theta,-\lambda;0+}^{-\delta} \mathbb{E}\, N^{\nu,\delta}
						(\cdot) \right) (t)
						\, \mathrm dt = \frac{\lambda^\delta}{s} + \lambda^\delta \int_0^\infty e^{-st} \mathbb{E} \,
						N^{\nu,\delta}(t) \, \mathrm dt \\
						& \Leftrightarrow \quad s^{\nu\delta+\theta} \int_0^\infty e^{-st} \int_0^t (t-y)^{\theta-1}
						E_{\nu,\theta}^{-\delta} \left[ -\lambda (t-y)^\nu \right]
						\mathbb{E}\,  N^{\nu,\delta} (y) \, \mathrm dy \, \mathrm dt
						=  \frac{\lambda^\delta}{s} + \lambda^\delta \int_0^\infty e^{-st} \mathbb{E} \,
						N^{\nu,\delta}(t) \, \mathrm dt \notag \\
						& \Leftrightarrow \quad s^{\nu\delta+\theta}
						\int_0^\infty \mathbb{E}\, N^{\nu,\delta}(y) \, \mathrm dy
						\int_{y}^\infty e^{-st} (t-y)^{\theta-1}
						E_{\nu,\theta}^{-\delta} \left[ -\lambda (t-y)^\nu \right]\, \mathrm dt
						= \frac{\lambda^\delta}{s} + \lambda^\delta \int_0^\infty e^{-st} \mathbb{E} \,
						N^{\nu,\delta}(t) \, \mathrm dt \notag \\
						& \Leftrightarrow \quad s^{\nu\delta+\theta}
						\int_0^\infty \mathbb{E}\, N^{\nu,\delta}(y) \, \mathrm dy
						\int_0^\infty e^{-s(z+y)} z^{\theta-1} E_{\nu,\theta}^{-\delta} (-\lambda z^\nu)\, \mathrm dz
						= \frac{\lambda^\delta}{s} + \lambda^\delta \int_0^\infty e^{-st} \mathbb{E} \,
						N^{\nu,\delta}(t) \, \mathrm dt \notag \\
						& \Leftrightarrow \quad s^{\nu\delta+\theta}
						\int_0^\infty e^{-sy} \mathbb{E}\, N^{\nu,\delta}(y) \, \mathrm dy
						\int_0^\infty e^{-sz} z^{\theta-1} E_{\nu,\theta}^{-\delta} \, \mathrm dz
						= \frac{\lambda^\delta}{s} + \lambda^\delta \int_0^\infty e^{-st} \mathbb{E} \,
						N^{\nu,\delta}(t) \, \mathrm dt \notag \\
						& \Leftrightarrow \quad s^{\nu\delta +\theta}
						\textsf L \left\{ \mathbb{E}\, N^{\nu,\delta} (t) \right\} (s) s^{-\theta}
						(1+\lambda s^{-\nu})^\delta = s^{-1} \lambda^\delta +\lambda^\delta
						\textsf{L} \left\{ \mathbb{E}\, N^{\nu,\delta} (t) \right\} (s) \notag \\
						& \Leftrightarrow \quad \textsf{L} \left\{ \mathbb{E}\, N^{\nu,\delta} (t) \right\} (s)
						= \frac{\lambda^\delta}{s \left[ (s^\nu+\lambda)^\delta -\lambda^\delta \right]}. \notag
					\end{align}
					Notice that the first step in \eqref{30mouse} is justified by the
					formula for the Laplace transform of the Riemann--Liouville fractional
					derivative and by applying the initial conditions.
					Before inverting the Laplace transform, it can be shown that
					\begin{align}
						\label{30acer}
						\textsf{L} \left\{ \mathbb{E} \, N^{\nu,\delta} (t) \right\} (s)
						& = \frac{\lambda^\delta}{s(s^\nu+\lambda)^\delta
						\left( 1-\frac{\lambda^\delta}{(s^\nu+\lambda)^\delta} \right)} \\
						& = \frac{\lambda^\delta}{s(s^\nu+\lambda)^\delta}
						\sum_{r=0}^\infty \left[ \frac{\lambda^\delta}{(s^\nu+\lambda)^\delta} \right]^r
						\notag \\
						& = \frac{1}{s} \sum_{r=0}^\infty \frac{\lambda^{\delta(r+1)}}{(s^\nu+\lambda)^{\delta(r+1)}}.
						\notag
					\end{align}
					Thus, the mean value is now easily found by inverting \eqref{30acer} term by term:
					\begin{align}
						\mathbb{E}\, N^{\nu,\delta} (t) & = \sum_{r=0}^\infty \int_0^t \lambda^{\delta(r+1)}
						y^{\nu\delta(r+1)-1} E_{\nu,\nu\delta(r+1)}^{\delta(r+1)} (-\lambda y^\nu) \, \mathrm dy \\
						& = \sum_{r=0}^\infty \lambda^{\delta(r+1)} t^{\nu\delta (r+1)}
					E_{\nu,\nu\delta(r+1)+1}^{\delta(r+1)}(-\lambda t^\nu). \notag
					\end{align}
				\end{proof}
		\end{Theorem}
		
		\begin{Remark}
			For $\delta = 1$, the mean value \eqref{30meanres} reduces to that of the pure fractional case
			\eqref{meanpure}. Indeed, 
			\begin{align}
				\mathbb{E} \, N^\nu(t) & = \sum_{r=0}^\infty \lambda^{r+1} t^{\nu(r+1)} E_{\nu,\nu(r+1)+1}^{r+1}
				(-\lambda t^\nu),
			\end{align}
			and passing now to the Laplace transform we get
			\begin{align}
				\textsf{L} \left\{ \mathbb{E}\, N^\nu (t) \right\} (s)
				= \sum_{r=0}^\infty \lambda^{r+1} s^{-\nu(r+1)} \left( 1+\lambda s^{-\nu} \right)^{-(r+1)}
				= \lambda / s^\nu,
			\end{align}
			which immediately leads to \eqref{meanpure}.  When $\delta=1$,   \eqref{30mean} reduces  to
			\begin{align}
				& \left[ \frac{\mathrm d^{\nu +\theta -k-1}}{\mathrm dt^{\nu +\theta -k-1}}
				\left( \bm{\mathrm E}_{\nu,\theta,-\lambda;0+}^{-1}
				\mathbb{E} \, N^\nu (\cdot) \right) (t) \right]_{t \to 0} = 0 \\
				& \Leftrightarrow \quad \left[ \frac{\mathrm d^{\nu +\theta -k-1}}{\mathrm dt^{\nu +\theta -k-1}}
				\int_0^t (t-y)^{\theta-1} E_{\nu,\theta}^{-\delta} [-\lambda(t-y)^\nu] \mathbb{E} \, N^\nu
				(y) \, \mathrm dy \right]_{t \to 0} = 0 \notag \\
				& \Leftrightarrow \quad \left[ \frac{\mathrm d^{\nu +\theta -k-1}}{\mathrm dt^{\nu +\theta -k-1}}
				\int_0^t (t-y)^{\theta-1} \left[ \frac{1}{\Gamma(\theta)}
				+\lambda (t-y)^\nu \frac{1}{\Gamma(1+\theta)} \right] \mathbb{E} \, N^\nu
				(y) \, \mathrm dy \right]_{t \to 0} = 0 \notag \\
				& \Leftrightarrow \quad\left[ \frac{\mathrm d^{\nu -k-1}}{\mathrm dt^{\nu -k-1}}
				\mathbb{E} \, N^\nu(t) + \lambda \mathbb{E} \, N^\nu(t) \right]_{t \to 0} = 0, \notag
			\end{align}
			for each $k= 0, \dots, n-1$, $n-1 \le \Re (\nu+\theta) < n$. By recalling $\mathbb{E}\,
			\mathcal{N}(t) = \sum_{r=0}^\infty r  p_r^\nu(t)$ and equation \eqref{40p} for $\delta=1$ we obtain
			\begin{align}
				\left[ \frac{\mathrm d^{\nu-1}}{\mathrm dt^{\nu-1}}
				\mathbb{E} \, N^\nu(t) \right]_{t \to 0} = 0,
			\end{align}
			where we considered only $k=0$  and therefore only one initial condition is used.
		\end{Remark}

		\subsection{Path simulation and parameter estimation}

			It is straightforward to generate a sample trajectory of generalization I by noting that the
			generalized Mittag--leffler random variable $T^{\nu,\delta}$ (see, e.g., \citet{pil})
			is a mixture of gamma densities, i.e.,
			\begin{equation}
				T^{\nu,\delta} \stackrel{d}{=} U^{1/\nu} V_\nu,
			\end{equation}
			where $U$ is gamma distributed with density function
			\begin{equation}
				f_U(u)=\frac{\lambda^\delta}{\Gamma (\delta)} u^{\delta-1} e^{-\lambda u}, \quad u>0,
			\end{equation}
			and $V_\nu$ is strictly positive-stable distributed with $\exp ( - s^\nu)$
			as the  Laplace transform  of the corresponding density function.
			Note that the $q$th fractional moment of the inter-event time can be easily shown as 
			\begin{equation}
				\mathbb{E}\left[ T^{\nu,\delta} \right]^q =\frac{\pi \Gamma (q/ \nu + 
				\delta)}{\lambda^{q / \nu} \Gamma ( q / \nu )\sin ( \pi q/ \nu) \Gamma (1- q )}, \qquad 0< q < \nu.
			\end{equation} 
			Typically, generating $T^{\nu,\delta}$ and adding one (corresponding to a single jump or event)
			each time gives a sample trajectory.

			Given $m$ jumps (corresponding to $m$ renewal times), we  propose  method-of-moments
			estimators for the  parameters  $\nu, \delta,$ and $\lambda$ to make the preceding
			generalization  usable in practice.  Getting the logarithm of $T^{\nu,\delta}$
			we have
			\begin{equation}
				T' \stackrel{d}{=} \frac{1}{\nu} U' +    V_\nu',
			\end{equation}
			where $T'=\ln (T^{\nu,\delta})$, $U'=\ln (U)$, and $V_\nu'=\ln (V_\nu)$.
			Following \citet{cuw} we get the estimating equations:
			\begin{equation}
				\mu_{T'} = \mathbb{E} \left( T' \right) =  \eta \bigg( \frac{1}{\nu} -1\bigg) +
				\frac{\psi ( \delta )-\ln (\lambda)}{\nu },
			\end{equation}
			\begin{equation}
				\sigma_{T'}^2= \frac{\pi^2}{6}\bigg( \frac{1}{\nu^2} -1\bigg) + \frac{1}{\nu^2} \psi^{(1)} (\delta),
			\end{equation}
			\begin{equation}
				\mu_3= \mathbb{E} \left( T'- \mu_{T'}  \right)^3 =
				\frac{ \psi^{(2)} (\delta) -2\left(\nu^3 -1\right)\zeta (3)}{\nu^3},
			\end{equation}
			$\eta \approx 0.57721$ is the Euler's constant, and   $\zeta (3)$ is the Riemann Zeta function
			evaluated at 3. Using the equations of the variance and the third central moment above,  we can solve for the
			estimates $\hat{\delta}$ and $\hat{\nu}$ using  $\hat{\mu}_3$ and $\hat{\sigma}_{T'}^2$. 
			Plugging $\hat{\nu}$ and $\hat{\delta}$ into the mean equation above, we obtain the estimate of $\lambda$ as
			\begin{equation}
				\hat{\lambda}=\exp \left(- \left[\hat{\nu} \left( \hat{\mu}_{T^{'}} -\eta (1/\hat{\nu}-1) \right)
				-\psi (\hat{\delta} ) \right] \right).
			 \end{equation}

			Furthermore, we tested the above procedure using the following estimate of the digamma function:
			\begin{equation}
				\psi(\tau) = \log (\tau) -  1/(2\tau)  - 1/(12\tau^2) + 1/(120\tau^4)  - 1/(252\tau^6) + O(1/\tau^8).
			\end{equation}
			We  then calculated the bias and the root-mean-square-error (RMSE) based on the 1000 generated
			data samples for different parameter values.  Table \ref{t2} in the appendix generally indicated positive
			results for the proposed method.

	\section{Generalization II}

		\label{00anothersec}	
		Recall that  a random variable $X$ is Mittag--Leffler-distributed
		with parameters $\lambda>0$ and $\nu \in (0,1]$ if it has probability density function
		\begin{align}
			\label{tonno}
			f_X(x) = \lambda x^{\nu-1} E_{\nu,\nu}(-\lambda x^\nu), \qquad x \in \mathbb{R}^+,
		\end{align}
		where
		\begin{align}
			E_{\alpha,\beta}(x) = \sum_{r=0}^\infty \frac{x^r}{\Gamma(\alpha r + \beta)}, \qquad x \in \mathbb{R},
		\end{align}
		is the Mittag--Leffler function. Note that $\text{Pr} \{ X>x \} = E_{\nu,1}(-\lambda x^\nu)$.
		
		Let  $Y=1/X$. Then the random variable $Y$ has the  inverse Mittag--Leffler distribution, that is, 
		\begin{align}
			\label{aa}
			\Pr \{ Y< y \} = \text{Pr} \{ 1/X < y \} = \text{Pr} \{ X > 1/y \} = E_{\nu,1}(-\lambda
			y^{-\nu}).
		\end{align}
		Hence, the corresponding probability density function is
		\begin{align}
			\label{a}
			\frac{\mathrm d}{\mathrm dy} \Pr \{ Y<y \} & = \frac{\lambda}{\nu} E_{\nu,\nu} (-\lambda y^{-\nu})
			\nu y^{-\nu-1} \\
			& = \lambda y^{-\nu-1} E_{\nu,\nu} (-\lambda y^{-\nu}), \qquad y \in \mathbb{R}^+. \notag
		\end{align}
		When $\nu=1$, formula \eqref{a} is the probability density function of an inverse
		exponential random variable, that is, 
		\begin{align}
			f(y) = \frac{\lambda}{y^2} e^{-\frac{\lambda}{y}}, \qquad y \in \mathbb{R}^+.
		\end{align}
		
\begin{figure}
			\centering
			\includegraphics[scale=.34]{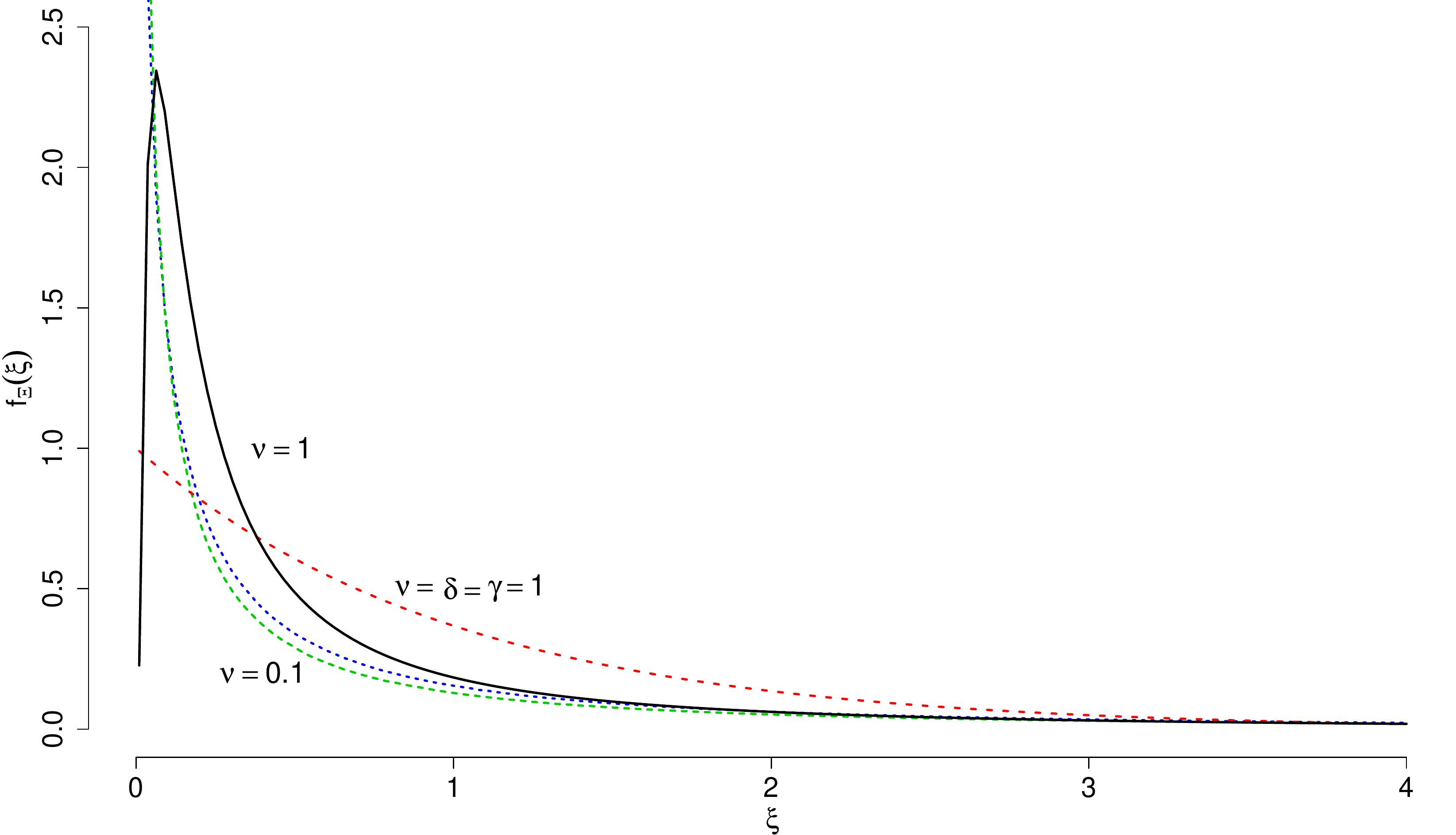} \\
			\includegraphics[scale=.34]{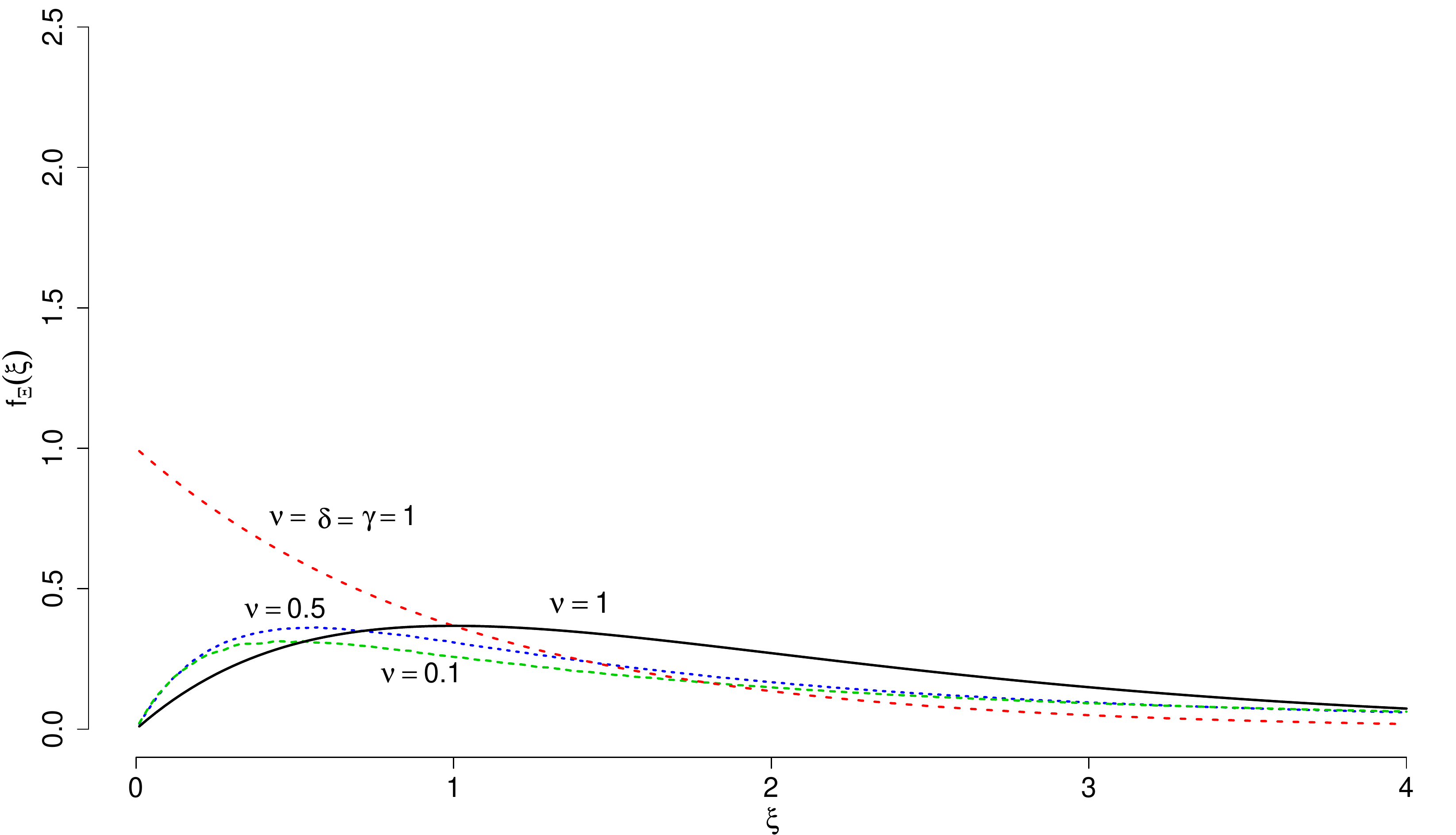} \\
			\includegraphics[scale=.34]{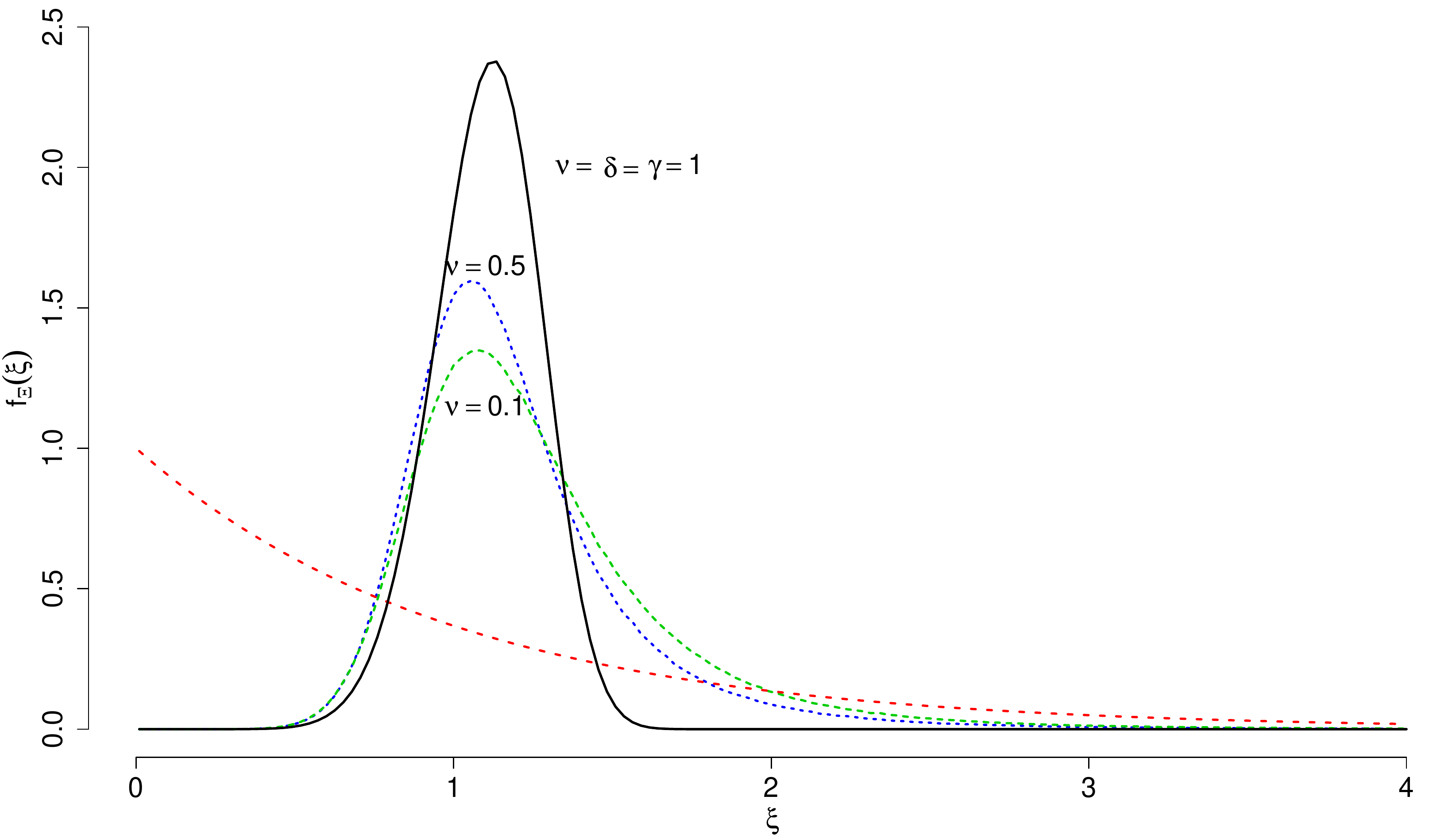}
			\caption{\label{afig3}The stretched-squashed Mittag--Leffler density plots (see \eqref{c}) for parameter values
			$(\nu, \delta, \lambda, \gamma )=((0.1, 0.5, 1), 2, 1, -0.5) $  (top),
			$(\nu, \delta, \lambda, \gamma )=((0.1, 0.5, 1), 2, 1, 1) $ (middle),
			and $(\nu, \delta, \lambda, \gamma )=((0.1, 0.5, 1), 2, 1, 5) $ (bottom).}
		\end{figure}

		We now give a single definition for both the Mittag--Leffler and the inverse
		Mittag--Leffler distributions. Note that the probability density function in equation \eqref{aa} can be written as
		\begin{align}
			\label{b}
			f_Y(y) = \lambda y^{\gamma -1} E_{\nu,\nu} (-\lambda y^\gamma), \qquad y \in \mathbb{R}^+, \: \nu
			\in (0,1],
		\end{align}
		where $\gamma = \pm \nu$.  By freeing the parameter $\gamma$
		in formula \eqref{b}, we arrive at the probability density
		\begin{align}
			\label{c}
			f_\Xi(\xi) = \frac{|\gamma|}{\nu} \lambda \xi^{\gamma-1} E_{\nu,\nu}(-\lambda \xi^\gamma),
			\qquad \xi \in \mathbb{R}^+, \: \nu\ \in (0,1], \: \gamma \in \mathbb{R} \backslash \{0\}
		\end{align}
		(see Figure \ref{afig3}).
		Observe  that $\int_0^\infty f_\Xi(\xi)\, \mathrm d \xi =1$  as
		\begin{align}
			\label{ciao}
			\frac{|\gamma|}{\nu} \lambda \int_0^\infty \xi^{\gamma -1} E_{\nu,\nu} (-\lambda \xi^\gamma)\,
			\mathrm d\xi
			& \overset{(\xi=z^{\nu/\gamma})}{=} \frac{|\gamma|}{\nu} \lambda
			\int_0^\infty \frac{\nu}{|\gamma|} z^{\nu-\nu/\gamma} E_{\nu,\nu}(-\lambda z^\nu)
			\, z^{\nu/\gamma-1} \mathrm dz \\
			& \quad = \int_0^\infty \lambda z^{\nu-1} E_{\nu,\nu} (-\lambda z^\nu)\,  \mathrm dz = 1. \notag
		\end{align}
		Note that in the second-to-last line of \eqref{ciao}, $\text{sgn} (\gamma)$ is used to stabilize the
		domain of the integral.
		
		 The Laplace transform $\mathbb{E} e^{-s \xi}$ can be shown as
		\begin{align}
			\label{acqua}
			\mathbb{E} e^{-s \xi} & = \lambda \frac{|\gamma|}{\nu} \int_0^\infty e^{-s \xi} \xi^{\gamma -1}
			E_{\nu,\nu} (- \lambda \xi^\gamma)\,  \mathrm d \xi \\
			& = \lambda  s^{-\gamma} \frac{|\gamma|}{\nu} \: {}_2 \psi_1 \left[ -\lambda s^{-\nu}
			\left|
			\begin{array}{l}
				(1,1),(\gamma,\nu) \\
				(\nu,\nu)
			\end{array}
			\right. \right] \notag \\
			& = \lambda s^{-\gamma} \frac{|\gamma|}{\nu} \sum_{r=0}^\infty \left( -\lambda s^{-\nu} \right)^r
			\frac{\Gamma (r \nu + \gamma)}{\Gamma(r \nu + \nu)},  \notag
		\end{align}
		where we use formula (2.2.22) of \citet{mathai}. When $\gamma = \nu \in (0,1]$ we obtain
		\begin{align}
			\mathbb{E} e^{-s \xi} = \lambda s^{-\nu} \sum_{r=0}^\infty \left( -\lambda s^{-\nu} \right)^r 
			= \frac{\lambda}{s^\nu + \lambda},
		\end{align}
		as in \citet{beg}, formula (4.15).
		
		From the above discussion it is clear that $\Xi = X^{\nu/\gamma}$, $\nu \in (0,1]$,
		$\gamma \in \mathbb{R}\backslash \{ 0 \}$,  where $X$ is the Mittag--Leffler distribution
		with probability density function \eqref{tonno}, therefore the waiting times $\Xi$ of a newly constructed
		renewal process are simple time-stretching or time-squashing of the original Mittag--Leffler
		distributed waiting times $X$.
		Below are density plots of $\Xi = X^{\nu/\gamma}$ where $X$ has the generalized Mittag--Leffler distribution
		given in \eqref{00fu}.

		Now, consider $m$ i.i.d.\ random inter-event imes $\Xi_1, \dots, \Xi_m$ of a
		counting process $\mathcal{N} (t)$, $t\ge 0$, which are distributed
		according to \eqref{c}.  		If $W_m = \Xi_1 + \dots + \Xi_m$ is the waiting time till the $m$th event, we have
		\begin{align}
			\mathbb{E} e^{-s W_m} = s^{-m \gamma} \left( \lambda \frac{|\gamma|}{\nu} \right)^m
			 {}_2 \psi_1^m \left[ -\lambda s^{-\nu}
			\left|
			\begin{array}{l}
				(1,1),(\gamma,\nu) \\
				(\nu,\nu)
			\end{array}
			\right. \right].
		\end{align}
		To determine the state probabilities $q_k^\nu(t) = \Pr \{ \mathcal{N}(t) = k \}$, $k \ge 0$,
		we write
		\begin{align}
			\int_0^\infty e^{-s t} q_k^\nu(t) \, \mathrm dt {} = &
			\int_0^\infty e^{-s t} \left( \Pr \{ W_k<t \} - \Pr \{ W_{k+1}<t \} \right) \mathrm dt \\
			= {} & s^{-k\gamma-1}\left( \lambda \frac{|\gamma|}{\nu} \right)^k
			{}_2 \psi_1^k \left[ -\lambda s^{-\nu} \left|
			\begin{array}{l}
				(1,1),(\gamma,\nu) \\
				(\nu,\nu)
			\end{array}
			\right. \right] \notag \\
			& - s^{-(k+1)\gamma-1}\left( \lambda \frac{|\gamma|}{\nu} \right)^{k+1} \! \! \!
			{}_2 \psi_1^{k+1} \left[ -\lambda s^{-\nu} \left|
			\begin{array}{l}
				(1,1),(\gamma,\nu) \\
				(\nu,\nu)
			\end{array}
			\right. \right]. \notag
		\end{align}

		\subsection{Path generation and parameter estimation}

			Simulating a sample path of generalization II directly follows from generalization I.
			Note that the $\Xi$'s can be generated using the algorithm of \citet{cuw}. It is also
			straighforward to show that the $q$th fractional moment of the random  inter-event time is
			\begin{equation}
				\mathbb{E}\mathcal{T}^q =\frac{\pi \Gamma (q \nu / \gamma + 1)}{\lambda^{q \nu / \gamma}
				\Gamma ( q / \nu )\sin ( \pi q/ \nu) \Gamma (1- q )}, \qquad 0< q < \nu.
			\end{equation}

			Given $m$ renewal times, we propose a formal procedure to estimate the parameters $\nu$, $\gamma$,
			and $\lambda$ of  generalization II. Let $\Xi' = \ln (\Xi)$ and $X' = \ln(X)$. Following \citet{cuw},
			we can deduce that
			\begin{equation}
				\label{e1}
				\mu_{\Xi'} = \mathbb{E} \left( \Xi' \right) =  \frac{\nu}{\gamma} \left( \frac{-\ln (\lambda) }{\nu}
				- \eta \right),
			\end{equation}
			\begin{equation}
				\label{e2}
				\sigma_{\Xi'}^2= \left( \frac{\nu}{\gamma} \right)^2\left[ \pi^2 \left( \frac{1}{3\nu^2} -\frac{1}{6}
				\right) \right],
			\end{equation}
			and
			\begin{equation} 
				\mu_3= \mathbb{E} \left( \Xi'- \mu_{\Xi'}  \right)^3 =  -2 \zeta (3) \left( \frac{\nu}{\gamma} \right)^3.
			\end{equation}
			Using the  estimating equations above, we can eliminate $\gamma$  and solve for $\nu$ by getting
			the $2/3$ root of the third central moment and dividing it by the variance. Thus, we obtain
			\begin{equation}
				\hat{\nu}= \sqrt{ \frac{c \pi^2}{3\left[ \left(2 \zeta (3) \right)^{2/3} + \frac{c \pi^2}{6}\right]} },
			\end{equation}
			where $c=( \hat{\mu}_3^{2/3} )/ \hat{\sigma}_{\Xi'}^2$.
			Substituting $\hat{\nu}$ to the variance equation \eqref{e2}, we get
			\begin{equation}
				\hat{\gamma}= \sqrt{ \left( \frac{\hat{\nu}}{\hat{\sigma}_{\Xi'}^2} \right)^2\left[ \pi^2 \left(
				\frac{1}{3 \hat{\nu}^2} -\frac{1}{6} \right) \right]  }.
			\end{equation}
			Finally, plugging $\hat{\nu}$ and $\hat{\gamma}$ into the mean equation \eqref{e1} above, we have
			\begin{equation}
				\hat{\lambda}= \exp \left[ - \left( \hat{\mu}_{\Xi'} \hat{\gamma} + \eta \hat{\nu}  \right)\right].
			\end{equation}

			We also tested the above explicit forms of the estimators  by calculating the bias and the
			root-mean-square-error (RMSE) based on the 1000 generated  data samples for different parameter
			and total jump size values.  Overall, Table \ref{t3} in the appendix showed favorable results for the
			proposed procedure.

	\section{Concluding remarks}

		\label{summary}
		We proposed two generalizations of the standard and the fractional Poisson processes through their
		renewal time distributions which naturally provided greater flexibility in modeling real-life renewal processes.	
		Statistical properties such as the state probabilities and process moments were derived.
		Algorithms to simulate
		trajectories and to estimate model parameters were also developed.
		Generally, tests provided additional merits to the
		proposed procedures.

		Although some work have already been done, there are still a few things that need to be pursued.
		For instance, the complete analysis of the counting process related to the renewal process
		that has stretched-squashed generalized
		Mittag--Leffler distributed waiting times would be a worthy pursuit.
		Also, the development of estimators using likelihood approaches would be of interest as well.

	\pagebreak

	\section{Appendix}
		\FloatBarrier

		\begin{table}[!ht]
			\centering
			\begin{tabular}{cc|ccc|ccc}
				\multicolumn{2}{c}{}  &  \multicolumn{3}{c}{Bias}  &  \multicolumn{3}{c}{RMSE} \\
				$\nu$ & $Est$ & $m=100$ & $1000$ & $10000$ & $m=100$ & $1000$ & $10000$ \\ \cline{1-8}
				\multirow{3}{*}{$0.5$}
				&  $\hat{\nu}$    & 0.018  & 0.001   & 0.000   &  0.184 & 0.049 & 0.013   \\
				&  $\hat{\delta}$ & 0.082  & 0.011   & 0.002   &  0.253 &  0.074 & 0.023     \\
				&  $\hat{\lambda}$ & 0.207  & 0.027   & 0.003   &  0.525 &  0.151 & 0.048     \\ \cline{1-8}
				\multirow{3}{*}{$0.6$}
				&  $\hat{\nu}$    & 0.026  & 0.003   & 0.000   &  0.210 & 0.068 & 0.016   \\
				&  $\hat{\delta}$ & 0.065  & 0.010   & 0.000   &  0.205 &  0.074 & 0.023     \\
				&  $\hat{\lambda}$ & 0.174  & 0.026   & 0.001   &  0.442 &  0.151 & 0.047     \\ \cline{1-8}
				\multirow{3}{*}{$0.7$}
				&  $\hat{\nu}$    & 0.024  & 0.001   & 0.000   &  0.254 & 0.056 & 0.018   \\
				&  $\hat{\delta}$ & 0.069  & 0.009   & 0.000   &  0.207 &  0.067 & 0.023     \\
				&  $\hat{\lambda}$ & 0.176  & 0.022   & 0.001   &  0.434 &  0.135 & 0.046     \\ \cline{1-8}
				\multirow{3}{*}{$0.8$}
				&  $\hat{\nu}$    & 0.005  & 0.003   & 0.000   &  0.264 & 0.077 & 0.020   \\
				&  $\hat{\delta}$ & 0.080  & 0.008   & 0.001   &  0.199 &  0.072 & 0.023     \\
				&  $\hat{\lambda}$ & 0.202  & 0.022   & 0.004   &  0.429 &  0.147 & 0.047     \\ \cline{1-8}
				\multirow{3}{*}{$0.95$}
				&  $\hat{\nu}$    & 0.022  & 0.001   & 0.000   &  0.349 & 0.079 & 0.021   \\
				&  $\hat{\delta}$ & 0.070  & 0.009   & 0.002   &  0.187 &  0.067 & 0.022     \\
				&  $\hat{\lambda}$ & 0.184  & 0.024   & 0.004   &  0.405 &  0.135 & 0.044     \\ \cline{1-8}
			\end{tabular}
			\caption{\emph{Parameter estimates for generalization I using different values of $\nu$,  $\delta=0.5$,
			and $\lambda=0.5$   for total jump sizes $m=100, 1000, 10000$.}}
			\label{t2}
		\end{table}

		\begin{table}[!ht]
			\centering
			\begin{tabular}{cc|ccc|ccc}
				\multicolumn{2}{c}{}  &  \multicolumn{3}{c}{Bias}  &  \multicolumn{3}{c}{RMSE} \\ 
				$\nu$ & $Est$ & $m=100$ & $1000$ & $10000$ & $m=100$ & $1000$ & $10000$ \\
				\cline{1-8}
				\multirow{3}{*}{$0.5$}
				&  $\hat{\nu}$     & 0.217   & 0.079   & -0.022  & 0.280 & 0.176   & 0.124 \\
				&  $\hat{\gamma}$ & -0.038   & -0.016   & 0.001  & 0.072 & 0.034   &  0.016  \\
				&  $\hat{\lambda}$ & -0.041   & -0.014   & 0.006  & 0.090 & 0.044   & 0.030   \\
				\cline{1-8}
				\multirow{3}{*}{$0.6$}
				&  $\hat{\nu}$    &  0.136   & -0.021   &0.000  & 0.222 & 0.153   &  0.100 \\
				&  $\hat{\gamma}$ & -0.026   & 0.018   & 0.051  & 0.067 & 0.032   &  0.016 \\
				&  $\hat{\lambda}$ & -0.017   & 0.005   & 0.001  & 0.085 & 0.041   &  0.023 \\
				\cline{1-8}
				\multirow{3}{*}{$0.7$}
				&  $\hat{\nu}$    & 0.052   &-0.022   &-0.008  & 0.188 & 0.149   & 0.060  \\
				&  $\hat{\gamma}$ &-0.015   & 0.001   & 0.001  & 0.071 & 0.034   & 0.013  \\
				&  $\hat{\lambda}$ &-0.003   & 0.007   & 0.002  & 0.086 & 0.040   & 0.014  \\
				\cline{1-8}
				\multirow{3}{*}{$0.8$}
				&  $\hat{\nu}$    &-0.013   &-0.030   &-0.003  & 0.172   & 0.129   & 0.034 \\
				&  $\hat{\gamma}$ & 0.004   & 0.006   & 0.001  & 0.073   & 0.035   &  0.011  \\
				&  $\hat{\lambda}$ & 0.007   & 0.008   & 0.000  & 0.086   & 0.035   & 0.009   \\
				\cline{1-8}
				\multirow{3}{*}{$0.95$}
				&  $\hat{\nu}$    &-0.043   &-0.004   & 0.000  & 0.131 & 0.044  & 0.013 \\
				&  $\hat{\gamma}$ & 0.021   & 0.002   & 0.000  & 0.083 & 0.028  & 0.008  \\
				&  $\hat{\lambda}$ & 0.014   & 0.002   & 0.000  & 0.073 & 0.022  & 0.007   \\
				\cline{1-8}
			\end{tabular}
			\caption{\emph{Parameter estimates for generalization II using different values of $\nu$,
			$\lambda=0.5$,  and  $\gamma=0.5$  for total jump sizes $m=100, 1000, 10000$.}}
			\label{t3}
		\end{table}

\pagebreak
	\bibliographystyle{unsrtnat}
	\bibliography{cahpol6}
	\nocite{*}		
		
\end{document}